\documentclass[a4paper]{amsart}

\usepackage{graphicx}
\usepackage{amsthm}
\usepackage{amssymb}
\usepackage{enumerate}
\usepackage{mathrsfs} 
\usepackage{hyperref}
\usepackage{subfig}
\usepackage{amsrefs}
\usepackage{mathtools}

\usepackage{oands}
\usepackage{hieroglf}

\theoremstyle{plain}%
\newtheorem{theorem}{Theorem}%  meant for continuous numbers
\newtheorem{question}[theorem]{Question}

\theoremstyle{remark}%
\newtheorem{remark}{Remark}%

\def\NN{\mathbb{N}}

\begin{document}

\title{On the Egyptian knotted rope}

\author{Rub\'en Vigara}

\address{IUMA, Universidad de
	Zaragoza}
\email{rvigara@unizar.es}

%    General info
\subjclass{01A16, 97F30, 97U60}

%\date{January 1, 2001 and, in revised form, June 22, 2001.}

\begin{abstract}
	The Egyptian knotted rope has been historically considered
	a measuring/geometric tool. A naive change in the way
	the rope is held to tighten it
	turns the rope into a primitive
	calculating tool --- a sort of primitive abacus.
	The arithmetic skills of the knotted rope have a double interest.
	From an historical point of view, they fit perfectly with some of the basic
	ancient Egyptian computing techniques.
	From an educational point of view, they enhance the potential
	interest of the knotted rope as a
	useful manipulative in mathematical education.
\end{abstract}

\keywords{knotted rope, arithmetic operations, Egyptian mathematics, mathematical manipulatives}


%\pacs[MSC Classification]{01A16, 97F30, 97U60}

\maketitle

\section{Introduction}\label{sec:intro}

The \emph{Egyptian knotted rope} or simply \emph{knotted rope}
or \emph{Egyptian rope}
is a rope or cord divided by knots into equal-length
segments.
It is usually presented as a geometric tool,
with a potential educational use for an
elementary introduction to right-angled triangles and the
Pythagorean theorem. Within the mathematical community,
the most widely held view regarding its use in ancient Egypt
is the construction of right triangles by
\emph{rope-stretching}, pulling from some knots based on 
Pythagorean triplets, although there is little evidence
of its actual use in this way~\cite{ImhausenMyths}.
However, there seems to be an academic consensus on the fact that 
knotted ropes --perhaps with their subdivisions marked with paint 
or another technique instead of knots-- were actually used by scribes in ancient Egypt,
at least for measuring purposes (see~\cite[p.154]{CRossiLibro} and references therein).
Indeed, the hieroglyphic sign for the number 100 is a coil of rope (\textpmhg{\HVi}).

A few years ago, while I was designing a geometry activity with the Egyptian rope
for K-3 classrooms, I realized that, besides its traditional use as a geometric tool,
the rope can also be used in a very simple way for basic arithmetic operations.
In particular, it can serve as a sort of ``primitive abacus'' for multiplications 
and divisions of small numbers
(such that the largest number involved in the computation is $\leq 100$, for instance).
Aside from their potential educational value, certain arithmetic operations
that are very easily carried out using the knotted rope strongly reminded me 
of some of the elementary computational techniques employed by scribes in ancient Egypt,
which I had previously heard about\footnote{The basic 
	principles of arithmetic in ancient Egypt were 
	explained to the author by A. M. Oller-Marc\'en.}.
These arithmetic properties of the Egyptian 
rope are so straightforward that one would expect 
them to have been mentioned before in the literature.
Unfortunately, I have found no mention of anything similar;
either I have been unable to find the correct reference,
or these ideas have simply been overlooked.
Since I believe these techniques should be widely known,
I shall describe them in this note.

I must emphasize that I am not an expert in either 
mathematics education or the history of mathematics;
whatever knowledge I have of these fields is elementary
and stems purely from personal interest.
I apologize in advance to experts in both fields for the
inaccuracies they will surely find in this
text from their specialist viewpoint. 
Whether the contents of this note have any
value for either of these disciplines
is for the experts to decide.

Although alternative ways of crafting the Egyptian rope ---using paint, 
fixed beads, buttons, or other types of markers instead of knots--- 
might be easier and even more accurate (to tie a bulky knot
in the exact desired position is not a trivial task),
in my opinion,
the use of knots produces a clean, elegant,
and visually attractive design using just one element---the rope.
The knot used in our real ropes
(Figures~\ref{fig:thick-cord} and~\ref{fig:myabacus})
is the Ashley knot
(``oysterman's knot'' in Ashley's book \cite[knot 526]{ashley}),
and the knot in the drawn ropes resembles
the double overhand knot~\cite[knot 516]{ashley}.

\section{Notation}

Our rope has a starting point, without a knot,
which can be identified with the number $0$,
and we will identify each knot with a natural
number that denotes its position when traveling along the rope
starting from $0$.
In this way, \emph{to identify} a given knot is to find the number that corresponds
to that knot, and \emph{to locate} the knot $n$ is the inverse operation:
to find the knot of the rope that corresponds to the number $n$.
A practical idea to simplify those tasks is to distinguish
the \emph{tens knots} (multiples of 10) by decorating them differently from the rest.
In our real-world pictures of Figures~\ref{fig:thick-cord} and~\ref{fig:myabacus}
we have also marked, with a different pattern, the \emph{fives knots} (odd multiples of 5).

At first, we consider no endpoint, so we can 
identify the rope with a \emph{flexible number line}, that may be bent
but may not be elongated or shortened. Given two knots $a,b\in\NN$, $a<b$, we denote
by $[a,b]$ the portion of the rope between $a$ and $b$, with both knots included.
In the same way, we denote by $[b,\infty)$, $b\in\NN$, the portion of the rope lying after $b$,
with $b$ included.

\subsection{Not-in-the-knot}

\begin{figure}
	\centering
	\hfill
	\subfloat[]{
		\includegraphics[width=0.25\textwidth]{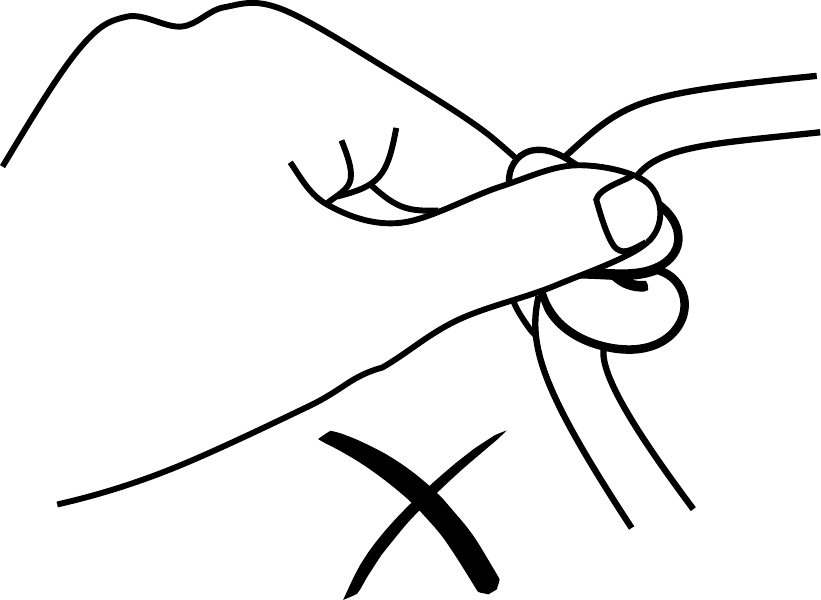}
		\label{fig:Not-in-the-knot-bad}}
	\hfill
	\subfloat[]{
		\includegraphics[width=0.3\textwidth]{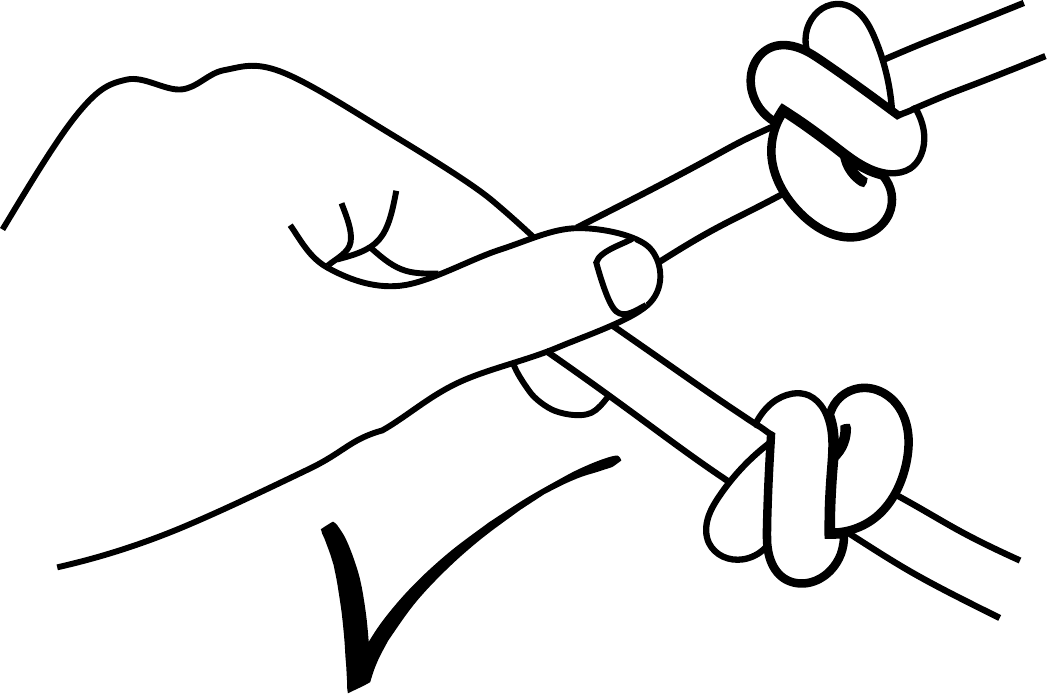}
		\label{fig:Not-in-the-knot-good}}
	\hfill\hfill
	\caption{Holding the rope: (a) at a knot, (b) at a midpoint between consecutive knots.}
	\label{fig:Not-in-the-knot}
\end{figure}

When using the rope as a length-measuring tool, the interval between
knots serves as the unit of measurement, and it is implicitly assumed
that the knots must lie at the endpoints of stretched segments.
Thus, there exists the widespread idea
that the rope-stretching operation should be done by pulling from
the knots (Fig.~\ref{fig:Not-in-the-knot-bad}).

The key to the arithmetic properties of the knotted rope
that we present here is to break with tradition at this point.
We will count knots, not intervals, so we hold
the rope at points \emph{between knots}, approximately halfway
between consecutive knots
(Fig.~\ref{fig:Not-in-the-knot-good}).

\subsection{Rope-stretching}
If $a\in\NN$ is a knot, to hold the rope \emph{before} or \emph{after} $a$ is to
hold the rope approximately at the midpoint between $a$ and $a-1$ or $a+1$, respectively.
If $a,b\in\mathbb{N}$, $a<b$ are two knots,
\emph{to stretch the segment} $[a,b]$ is the operation that consists in holding the
rope before $a$ and after $b$, pulling in opposite directions until the rope between
the two holding points becomes tight.

\section{Arithmetic with the knotted rope}

Obviously, addition and subtraction correspond 
to counting on the rope forward and backward, respectively.
The operation that interests us most is the following.
\subsection{Cloning segments (duplication)}
Segment stretching allows immediate visual comparisons between numbers/segments.
Thus, we can \emph{clone} segments in quite a simple way.
If we have stretched the segment $[1,n]$, in order to obtain the segment $[n+1,2n]$,
or to locate the knot $2n$, it suffices to
pull the next section $[n+1,\infty)$ of the rope
in the opposite direction of $[1,n]$ until we
obtain a stretched segment equal to the original one
(Figure~\ref{fig:doubling}).
\begin{remark}
	Once the original segment is stretched, 
	the cloning operation is \emph{computationally free}
	in the sense that it does not require any arithmetic
	or counting operation to be carried out.
\end{remark}

\subsection{Multiplication}
\begin{figure}
	\centering
	\subfloat[doubling/halving]{
		\includegraphics[width=0.9\textwidth]{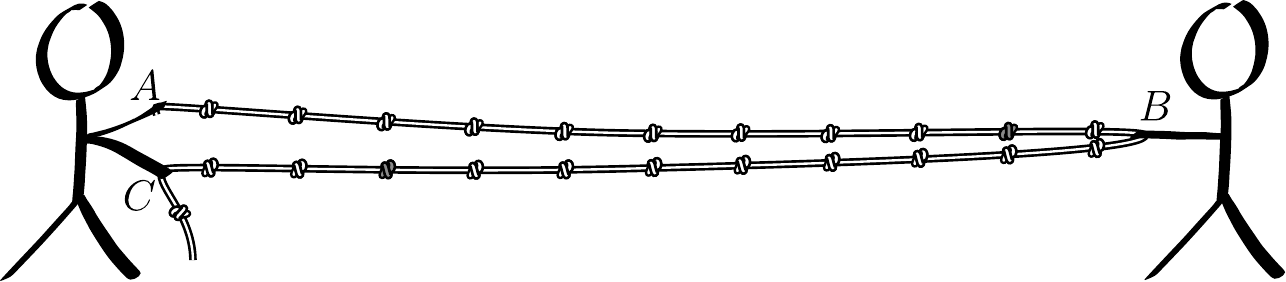}
		\label{fig:doubling}}
	\\
	\subfloat[tripling/trisecting... taking $2/3$?]{
		\includegraphics[width=0.9\textwidth]{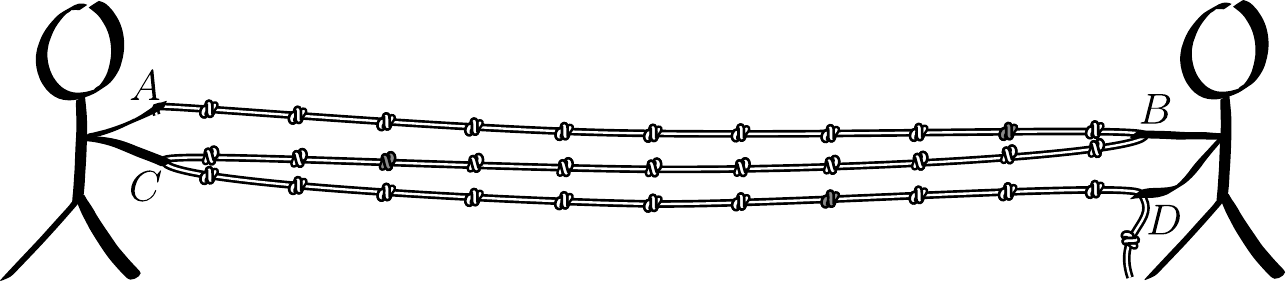}
		\label{fig:tripling}}
	\caption{Basic arithmetic operations with the knotted rope.}
	\label{fig:dummies}
\end{figure}

The cloning operation can be performed
over and over again in order to obtain the subsequent multiples of the original segment,
leading to a \emph{zigzag array} of cloned stretched segments
(Figure~\ref{fig:tripling} and~\ref{fig:abacus}).
In our figures, the starting point of the rope is always at the
top left corner of the zigzag array.
Thus, to compute $m \times n$: (i) stretch the segment $[1,n]$;
(ii) clone it $m-1$ times; (iii) identify the last knot of the
last segment of the zigzag array so obtained.

\subsection{Division}
Integer division can also take advantage of the cloning operation.
To compute $m \div n$, $m>n$, in order to obtain the quotient $q$
and the remainder $r$,
a tentative algorithm using the rope
would be as follows: (i) locate the dividend knot $m$ 
and mark it in some way (with a clothespin or something else);
(ii) locate the divisor knot $n$ and stretch $[1,n]$;
(iii) clone $[1,n]$ and repeat cloning while $m$ is
outside the zigzag array; (iv) stop cloning after the last cloning
before $m$ gets into the array; (v) $q$
is the number of stretched segments in the zigzag array, and $r$
is the number of knots outside the array, up to (and including) $m$.

\subsection{Division II. Halving and trisecting}
In the previous paragraph the quotient
$m \div n$ is the result of a quotitive division
(``the number of groups of $n$ elements that can be formed 
from a set of $m$ elements'').
The equal-sharing or partitive division $m \div n$ (``size of each part
after splitting $m$ into $n$ equal parts'') can be performed
so easily with the knotted rope in the simplest cases
$n=2,3$ (Figure~\ref{fig:dummies}): \emph{halving} and \emph{trisecting}. 

To halve $n$ (Figure~\ref{fig:doubling}): (i) hold the rope 
before $1$ and after $n$ (points $A$ and $C$, respectively)
and keep these holding points close together;
(ii) pull from a sliding pivot point $B$ between $1$ and $n$
away from the holding points, until the rope between $A$ and $C$ becomes tight;
(iii) identify the last knot in the first stretched segment.

To trisect $n$ (Figure~\ref{fig:tripling}): (i) hold the rope 
before $1$ and after $n$ (points $A$ and $D$, respectively), 
and insert two sliding pivot points $B$ and $C$ between the holding points;
(ii) keep together $A$ with $C$, and $B$ with $D$, and pull $B,D$ away from
$A,C$ until the rope between $A$ and $D$ becomes tight as in
Figure~\ref{fig:tripling};
(iii) identify the last knot in the first stretched segment.

Halving and trisecting procedures
can be easily adapted for the cases when
$n/2$ or $n/3$ respectively is not integer, whether we are interested in
the integer division or in the exact division result.

\begin{remark}\label{remark:cloning-without-stretching}
	Doubling/halving can be carried out by a single person ``locally'', i. e.,
	without moving from their position and without fully stretching a whole segment.
\end{remark}

For example, to clone the segment $[a,b]$:
(i)	put a clothespin before $a$ ($a$\emph{--clothespin});
(ii) fold the rope halfway between $b$ and $b+1$ and put another clothespin at the bend ($b$\emph{--clothespin});
(iii) hold both strands of the rope with one hand ---let's say your left--- near the 
$b$--clothespin, but loosely,
so that the rope can slide through your hand;
(iv) with the free (right) hand, pull both strands of the rope, starting from the bend,
so that the knots pass through your left hand in pairs, one from each strand;
(v) repeat pulling until the $a$--clothespin reaches your left hand;
(vi) put a last clothespin ($c$\emph{--clothespin}) in the 
strand $[b+1,\infty)$ in the parallel position
of the $a$--clothespin. The segment between the $b$-- and
$c$--clothespins is a clone of $[a,b]$.
An equivalent procedure can be used for halving.

\section{Related questions}
As the reader will have noted, the mathematics discussed
in this note are highly basic. So basic that,
as we have mentioned before, they
should have been previously known by the community.
Therefore, the first question we raise is the following:
\begin{question}	
	Why have we not seen these constructions before?
\end{question}

\subsection{Ancient Egyptian arithmetic}

Scribes in ancient Egypt
performed multiplication and division through 
different combinations of a few basic operations.
Among them, some of the most commonly used were
doubling/halving, decupling and taking $2/3$
(see~\cite{gillings,imhausenLibro}).

Decupling could arise naturally from the use of a decimal number system
but, as far as I know, it is not clear why the other operations
were so common.

As we have seen, the doubling/halving operation in particular
(see Remark~\ref{remark:cloning-without-stretching}), at least for
numbers up to 100, can be performed with relative ease with the knotted rope.

What about taking $2/3$?
\begin{remark}\label{remark:dummy-two-thirds}
	Assume a trisecting operation as in Figure~\ref{fig:tripling}.
	If,	for any reason (right dummy holds a lower rank,
	is innumerate,... or simply does not exist at all!), the dummy on the left
	is in charge of noting down the result of the operation,
	they will have direct access to $\frac23 n$ 
	instead of $\frac13 n$.
\end{remark}
Thus, at least theoretically, the ``taking--$2/3$'' operation might also 
have arisen through the use of a rope as a computing tool.
\begin{question}
	Was the knotted rope used in ancient Egypt, perhaps
	in an early stage of their mathematical development, as an
	arithmetical tool in a similar manner to that
	described in this note?
\end{question}

\subsection{Mathematical learning: the knotted rope abacus}

\begin{figure}
	\centering	
	\subfloat[]{
		\includegraphics[width=0.85\textwidth]{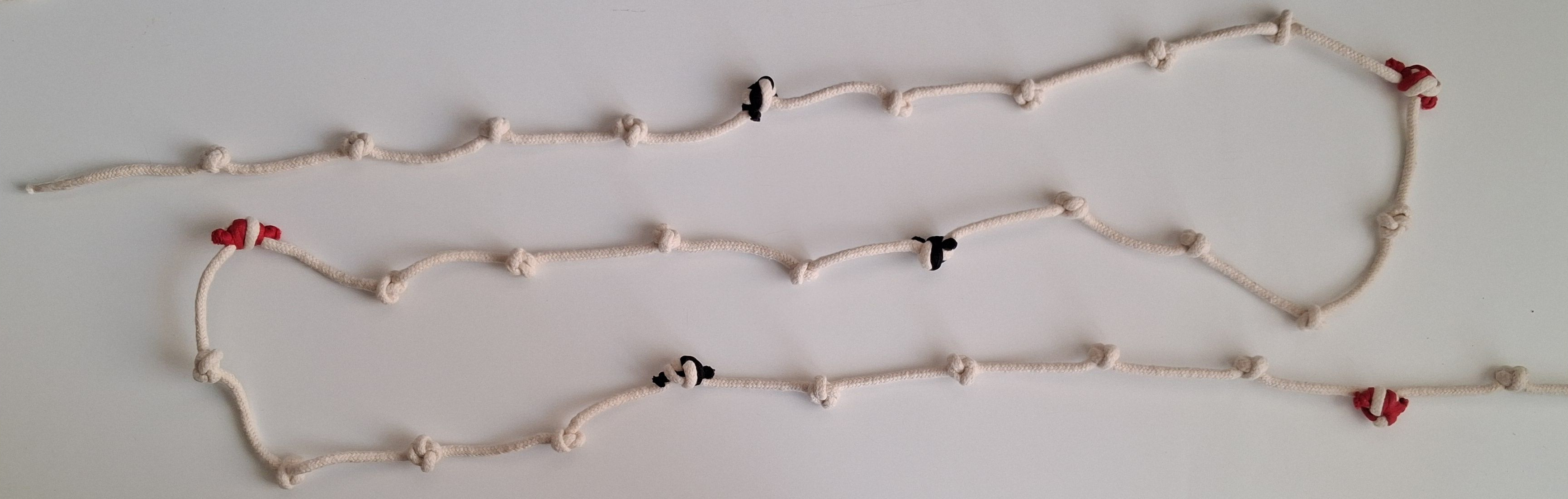}
		\label{fig:thick-cord}}\\
	\subfloat[]{
		\includegraphics[width=0.45\textwidth]{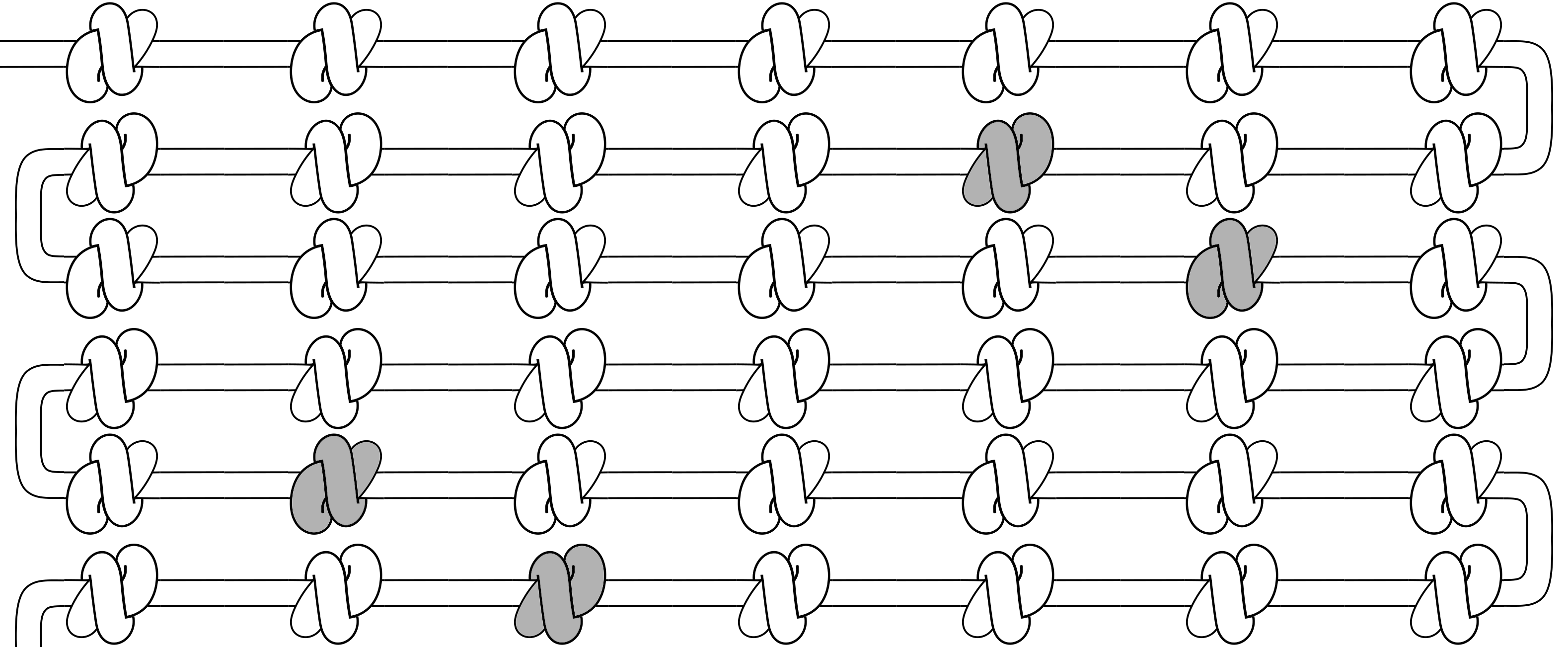}
		\label{fig:virtual6x7}}\hfill
	\subfloat[]{
		\includegraphics[width=0.45\textwidth]{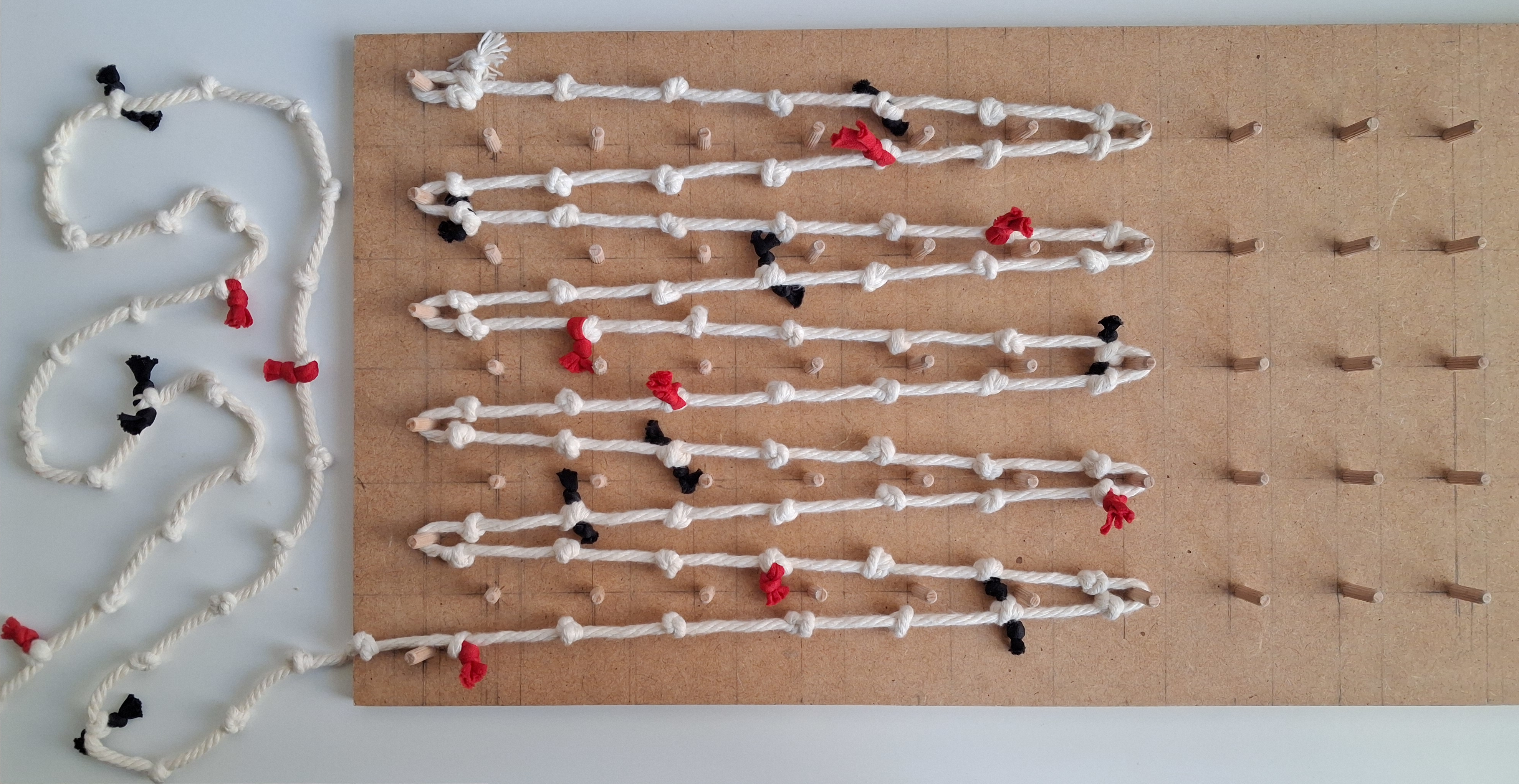}
		\label{fig:myabacus}}
	\caption{(a) A thick knotted rope. (b) ``Virtual'' knotted rope abacus for computing $6\times 7$. (c) A real prototype displaying
		the table of multiples of 7.}
	\label{fig:abacus}
\end{figure}
See~\cite{Hurst2020} for a quick review of the literature
about the use of manipulatives in mathematical instruction,
with special attention on multiplicative thinking.

It seems reasonable to think
that the knotted rope could be introduced in school as a manipulative 
from the earliest stages in such a way that the four principles
for maximizing the effectiveness of manipulatives
proposed in~\cite{Laski2015} are satisfied.

A thick knotted rope with 20--30 bulky knots dividing
the rope into intervals of
50~cm, for example\footnote{It is assumed that the ropes in ancient Egypt
	were divided into 1-cubit ($\approx$ 52.5~cm) 
	intervals~\cite{imhausenLibro,CRossiLibro}.},
with tens or tens and fives knots decorated with specific 
patterns (Figure~\ref{fig:thick-cord}),
could be used by young students in a collaborative way, as suggested
in Figure~\ref{fig:dummies}, for measuring, geometric and arithmetic tasks,
providing a first contact with the number line. 

For older students, a thinner cord with one hundred knots with a
separation of 10~cm, for example, combined with
a device allowing the simple construction
of zigzag arrays of stretched segments of up to ten knots,
could serve as a sort of abacus to assist in learning multiplication tables.
We show a virtual version of how such \emph{knotted rope abacus}
should work in
Figure~\ref{fig:virtual6x7}, and a real prototype in
Figure~\ref{fig:myabacus}.
It is fairly easy to visualize the product $6\times 7$,
for example, with Cuisenaire rods, a 100-bead abacus, arrays or other
well-known manipulatives,
but these tools do not eliminate the need for 
repeated addition or counting to determine the result of the operation.
A zigzag array as the one of Figure~\ref{fig:virtual6x7}
or the top 6 stretched segments of Figure~\ref{fig:myabacus},
combines a nice visualization with computational efficiency,
as it leads to an easy determination of the result of the operation.

If the markings or decorations that distinguish 
the tens (or tens and fives) knots incorporate tactile or even 
auditory differences (such as small bells) compared to 
the other knots, the rope could even be used 
by students with visual impairments
(Remark~\ref{remark:cloning-without-stretching}).

Thus, our last question is:
\begin{question}
	Can this manipulative be used to help students develop fluency in multiplicative thinking
	or other areas of mathematical learning?
\end{question}

%\section*{Declarations}
%
%\begin{itemize}
%	\item Funding. The author has been supported by the Spanish Research project PID2024-156032NB-I00,
%	funded by MICIU/AEI/10.13039/501100011033 and FEDER, UE,
%	and by European Regional Development Fund and Diputaci\'on General de Arag\'on (ER22\_23R).
%	\item Competing interests. The author has no competing interests to declare that
%	are relevant to the content of this article.
%	\item Ethics approval and consent to participate. `Not applicable'
%	\item Consent for publication. `Not applicable'
%	\item Data availability. `Not applicable'
%	\item Materials availability. `Not applicable'
%	\item Code availability. `Not applicable'
%	\item Author contribution. `Not applicable'
%\end{itemize}
%
%\noindent
%If any of the sections are not relevant to your manuscript, please include the heading and write `Not applicable' for that section. 

\bibliographystyle{amsplain}
\bibliography{bibliography2}% common bib file

\end{document}